\documentclass[12pt, reqno]{amsart}
\usepackage [latin1]{inputenc}
\usepackage{amscd}
\usepackage{epsfig}
\usepackage{amssymb}
\usepackage{amsmath}
\usepackage{amsthm}
\usepackage[T1]{fontenc}
\usepackage{ae,aecompl}

\usepackage {amsfonts} 
\usepackage {latexsym} 
\usepackage {exscale} 
\usepackage{amsmath}

\newcommand{\R} {\ensuremath{\mathbb{R}}}

\newcommand{\C} {\ensuremath{\mathbb{C}}}
\newcommand{\Z} {\ensuremath{\mathbb{Z}}}

\newcommand{\OO}{\mathcal{O}}
\renewcommand{\o}[1]{\overline{#1}}

\newcommand{\dq}{\overline{\partial}}
\newcommand{\wt}[1]{\widetilde{#1}}
\DeclareMathOperator{\Reg}{Reg}

\newtheorem {satz} {Satz} [section]
\newtheorem {lem} [satz] {Lemma}

\newtheorem {defn} [satz] {Definition}

\newtheorem {thm} [satz] {Theorem}

\DeclareMathOperator{\dist}{dist}
\DeclareMathOperator{\supp}{supp}

\renewcommand{\theta}{\vartheta}

\begin{document}

\title[$\dq$ on homogeneous varieties] 
{The $\dq$-equation on homogeneous varieties with an isolated singularity}

\author{J. Ruppenthal}

\address{Department of Mathematics, University of Wuppertal, Gau{\ss}str. 20, 42119 Wuppertal, Germany.}

\email{ruppenthal@uni-wuppertal.de}

\date{July 23, 2008}


\begin{abstract}
Let $X$ be a regular irreducible variety in $\C\mathbb{P}^{n-1}$, $Y$ the associated homogeneous variety
in $\C^n$, and $N$ the restriction of the universal bundle of $\C\mathbb{P}^{n-1}$ to $X$.
In the present paper, we compute the obstructions to solving the $\dq$-equation in the $L^p$-sense on $Y$
for $1\leq p\leq \infty$ in terms of cohomology groups $H^q(X,\OO(N^\mu))$.
That allows to identify obstructions explicitly if $X$ is specified more precisely,
for example if it is equivalent to $\C\mathbb{P}^1$ or an elliptic curve.
\end{abstract}


\maketitle


\section{Introduction}

One strategy to study the $\dq$-equation
on singular complex spaces is to use Hironaka's resolution of singularities in order to pull-back
the $\dq$-equation to a regular setting, where it is treatable much easier. 
See \cite{AHL}, \cite{BiMi} or \cite{Ha} for 
detailed information about resolution of singularities. That strategy has been pursued already
in \cite{FOV1} and \cite{Rp4}, where it leads to more or less imprecise results.
But the method seems to be quite promising for further investigations,
because it can be improved considerably.
We were able to do that in this paper for homogeneous varieties
with an isolated singularity, where the desingularization is obtained by a single blow up.
We believe that one should draw special attention to this strategy,
because there are some analogies to the case of complex projective varieties,
where we have an intimate connection between 
the $L^2$-cohomology of the regular part of the variety
and the $L^2$-cohomology of resolutions (see \cite{PaSt1}).

\vspace{2mm}
For a complex projective variety $Z\subset\C\mathbb{P}^n$,
the Cheeger-Goresky-MacPherson conjecture (see \cite{CGM}) states that
the $L^2$-deRham cohomology  $H^*_{(2)}(Z^*)$
of the regular part of the variety $Z^*:=\mbox{Reg } Z$
with respect to the (incomplete) restriction of the Fubini-Study metric
is naturally isomorphic to the intersection cohomology of middle perversity $IH^*(Z)$
(which in turn is isomorphic to the cohomology of a small resolution of singularities).
Ohsawa proved this conjecture under the extra assumption that the variety has only isolated singularities in \cite{Oh},
while it is still open for higher-dimensional singular sets.
The early interest in the conjecture of Cheeger, Goresky and MacPherson was motivated in large
parts by the hope that one could then use the natural isomorphism and a Hodge decomposition
for $H^k_{(2)}(Z^*)$ to put a pure Hodge structure on the intersection cohomology of $Z$ (cf. \cite{CGM}). That
was in fact done by Pardon and Stern in the case of isolated singularities (see \cite{PaSt2}).
Their work includes the computation of the $L^2$-Dolbeault cohomology groups $H^{p,q}_{(2)}(Z^*)$
in terms of cohomology groups of a resolution of singularities (see \cite{PaSt1}; also for further references).\\

Let us now direct our attention to the case of singular Stein spaces.
Though one would expect similar relations in this (local) situation,
no such representation of the $L^2$-Dolbeault cohomology is known.
The best results include quite rough lower and upper bounds on the dimension
of some of the $L^2$-Dolbeault cohomology groups (see \cite{DFV}, \cite{Fo}, \cite{FOV1}, \cite{FOV2}, \cite{OvVa} or \cite{Rp4}).
The origin of the present work is the attempt to compute the $L^2$-Dolbeault cohomology groups
in the spirit of the work of Cheeger-Goresky-MacPherson, Ohsawa, Pardon-Stern and others
in terms of certain cohomology groups on a resolution of singularities.
But, in the absence of compactness, most of their arguments do not carry over to the local situation
and one has to develop some new strategies.

\vspace{2mm}
One such tool which could be helpful for studying the $\dq$-equation (even locally) on singular complex spaces
is a Dolbeault complex with weights according to normal crossings developed in \cite{Rp6}.
A short review of this construction is contained in section \ref{sec:sufficient} of this paper,
the main result is the exactness of the complex cited here as Theorem \ref{thm:main}.
Weights according to normal crossings are a natural choice because
one can achieve that the exceptional set of a desingularization consists of normal crossings only,
and the deformation of a metric under desingularization produces singularities along the exceptional set
which have to be taken into account when we treat the $\dq$-equation (cf. the introduction to \cite{Rp6}).

\vspace{2mm}
Another interesting tool that we use in this paper is an integration along the fibers
of the normal bundle of the exceptional set of a desingularization.
This idea has been already used by E.\ S.\ Zeron and the author in \cite{RuZe} to construct an explicit
$\dq$-integration formula on weighted homogeneous varieties.
The method is described in section \ref{sec:fibers}.
A crucial point about both these tools is that they depend on integral formulas.
So, they allow to drop the
the restriction to $L^2$-spaces given by the well-known Hilbert space methods.

\vspace{2mm}
In view of the large difficulties in computing the $L^2$-cohomology explicitly, 
it seems reasonable
to gain a broader view and better understanding by also considering $L^p$-Dolbeault cohomology groups
for arbitrary $1\leq p\leq \infty$. 
Besides the $L^2$-results mentioned above,
only the $L^\infty$-case has been addressed in a number of publications:
\cite{AcZe1}, \cite{AcZe2}, \cite{FoGa}, \cite{Rp2}, \cite{Rp4}, \cite{RuZe}, \cite{SoZe}.
These papers treat H{\"o}lder regularity of the $\dq$-equation provided the right-hand side
of the equation is bounded. Clearly, this implies the solution of the Cauchy-Riemann equations in the $L^\infty$-sense.
In view of those results, the present paper is an attempt to embed the $L^2$ and $L^\infty$-case into
the broader spectrum of an $L^p$-theory.
In fact, by use of the Dolbeault complex with weights and the integration along the fibers of the normal bundle,
is is possible to compute the $L^p$-Dolbeault cohomology groups on a homogeneous variety
with an isolated singularity $Y$ for all $p$ such that $2d/p\notin \Z$ (where $d=\dim Y$) and for $p=1$
(see Theorem \ref{thm:sufficient} and Theorem \ref{thm:necessary} below).
This does not solve the $L^2$-problem but gives a quite precise idea what to expect for the $L^2$-groups.\\


\vspace{2mm}
We will now describe the results of this paper in detail.
Let $X$ be a regular irreducible variety in $\C\mathbb{P}^{n-1}$, and $Y$ the associated homogeneous variety
in $\C^n$ which has an isolated singularity at the origin. We denote by $N$ the restriction of the universal bundle
on $\C\mathbb{P}^{n-1}$ to $X$. Let $d=\dim Y$.

\vspace{2mm}
The regular complex manifold $Y^*:=Y\setminus \{0\}=\Reg Y$ carries a hermitian structure induced
by restriction of the euclidian metric of the ambient space $\C^n$.
Let $|\cdot|_Y$ and $dV_Y$ be the resulting metric and volume form on $Y^*$.
Now, if $U\subset Y^*$ is an open set, and $\omega$ a measurable $(0,q)$-form on $U$,
we set
\begin{eqnarray*}
\|\omega\|_{L^p_{0,q}(U)}^p &:=& \int_U |\omega|_Y^p dV_Y\ , \mbox{ for } 1\leq p<\infty,\\
\|\omega\|_{L^\infty_{0,q}(U)} &:=& \mbox{ess} \sup_{\substack{z\in U}} |\omega|_Y(z).
\end{eqnarray*}

We are interested in the following cohomology groups,
where the $\dq$-equation has to be interpreted in the sense of distributions
(throughout this paper).
Due to the incompleteness of the metric, different extensions of the $\dq$-operator on smooth forms
lead to different cohomology groups. For $U\subset \Reg Y$ open, let
\begin{eqnarray*}
H^q_{(p)}(U,\OO) := \frac{\{\omega\in L^p_{0,q}(U):\dq\omega=0\}}{\{\omega\in L^p_{0,q}(U):\exists f\in L^p_{0,q-1}(U): \dq f=\omega\}}.
\end{eqnarray*}

We will show (giving sufficient conditions for $L^p$-solvability of the $\dq$-equation):

\begin{thm}\label{thm:sufficient}
Let $X$, $Y$ and $N$ as above, $D\subset\subset Y$ strongly pseudoconvex such that $0\in D$, $D^*=D\setminus\{0\}$,
and $1\leq p\leq \infty$, $1\leq q\leq d=\dim Y$. Set
$$a(p,q,d) :=
\left\{\begin{array}{ll}
\max\{k\in\Z: k< 1+q- 2d/p\} &, p \neq 1,\\
\max\{k\in\Z: k\leq 1+ q -2d/p\} &, p=1.
\end{array}\right.$$
Then there exists an injective homomorphism
\begin{eqnarray}\label{eq:sufficient}
H^q_{(p)}(D^*,\OO) \hookrightarrow \bigoplus_{\mu\geq a(p,q,d)} H^q(X,\OO(N^{-\mu})).
\end{eqnarray}
\end{thm}

The right hand side in \eqref{eq:sufficient} is finite-dimensional because $N$ is a negative holomorphic line bundle.
Necessary conditions are determined by:

\begin{thm}\label{thm:necessary}
Let $X$, $Y$ and $N$ as above, and let $D\subset\subset Y$ be an open set such that $0\in D$, $D^*=D\setminus\{0\}$,
and $1\leq p\leq \infty$, $1\leq q\leq d=\dim Y$. Set
$$c(p,q,d) :=
\max\{k\in\Z: k\leq 1+q- 2d/p\}.$$
Then there exists an injective homomorphism
\begin{eqnarray}\label{eq:necessary}
 \bigoplus_{\mu\geq c(p,q,d)} H^q(X,\OO(N^{-\mu})) \hookrightarrow H^q_{(p)}(D^*,\OO).
\end{eqnarray}
\end{thm}


Note that sufficient and necessary conditions coincide if $2d/p\notin\Z$ or $p=1$,
and that $c(p,q,d)=a(p,q,d)+1$ in all other cases. So, there remains a little uncertainness
about the contribution of $H^q(X,\OO(N^{-a}))$, for example if $p=2$.\\


The proof of Theorem \ref{thm:sufficient} and Theorem \ref{thm:necessary}
depends heavily on an embedded desingularization of $Y\subset \C^n$, which is in our situation simply given by
a single blow-up of the origin in $\C^n$.
We will study the behavior of $L^p$-norms under this resolution of singularities in the next section,
while we will present the first part of the proof of Theorem \ref{thm:sufficient} in section \ref{sec:sufficient}.
The main tool here is a Dolbeault complex with weights according to normal crossings 
that was constructed in \cite{Rp6}.
The second part of the proof is settled by
another important tool of our work, namely an integration along the fibers of the holomorphic line bundle $N$,
which we will develop in section \ref{sec:fibers}.
This idea has been already used by E.\ S.\ Zeron and the author in \cite{RuZe} to construct an explicit
$\dq$-integration formula on weighted homogeneous varieties.
In section \ref{sec:fibers}, we obtain as a byproduct:

\begin{thm}\label{thm:integration}
Let $X$ and  $Y$ be as above, $D\subset\subset Y$ an open subset, $D^*=D\setminus\{0\}$,
and $1\leq p\leq \infty$, $1\leq q\leq \dim Y$. Let
$\omega \in L^p_{0,q}(D^*)\cap \ker \dq$
with compact support in $D$. Then there exists
$\eta\in L^p_{0,q-1}(D^*)$
such that $\dq\eta=\omega$.
\end{thm}

Using Theorem \ref{thm:integration} in case $q=1$ and Hartogs' Extension Theorem on normal Stein spaces
with isolated singularities, it is easy to deduce vanishing of the first cohomology with compact support (see section \ref{sec:fibers}):

\begin{thm}\label{thm:compact}
Let $X$ and  $Y$ be as above, $D\subset\subset Y$ an open subset, $D^*=D\setminus\{0\}$,
and $1\leq p\leq \infty$. Then:
$$H^1_{(p),cpt}(D,\OO):=\frac{\{\omega\in L^p_{0,1}(D^*): \dq\omega=0,\ \supp\omega\subset\subset D\}}
{\{\ \omega\in L^p_{0,1}(D^*): \exists f\in L^p(D^*): \dq f=\omega,\ \supp f\subset\subset D\}}=0.$$
\end{thm}

\vspace{1mm}
We will then prove Theorem \ref{thm:necessary} in section \ref{sec:necessary},
and discuss some examples and applications in the last section \ref{sec:examples}.
Let us mention a few of them at this point. Let $X$, $Y$ and $N$ be as above, and $D$ a strongly pseudoconvex neighborhood of the origin in $Y$, $D^*=D\setminus\{0\}$.
If, for example, a group $H^q_{(p)}(D^*,\OO)$ is vanishing, then it follows by standard techniques
that we can construct a bounded $L^p$-solution operator for the $\dq$-equation in degree $(0,q)$ on $D^*$ (see Theorem \ref{thm:bounded}).

\vspace{2mm}
When we restrict our attention to the case $\dim Y=2$, $X$ is a compact Riemann surface,
and that allows to compute the groups $H^1(X,\OO(N^{-\mu}))$ by the Theorem of Riemann-Roch.
We will do that for $X\cong \C\mathbb{P}^1$ or $X$ an elliptic curve, and deduce some consequences for $L^p$-solvability of
the $\dq$-equation on $Y$.

\vspace{2mm}
Combining an Extension Theorem for cohomology classes on complex spaces of Scheja (Theorem \ref{thm:scheja1}) with our
integration along the fibers, we deduce that
$$H^q_{(p)}(D^*,\OO)=0$$
for $1\leq q \leq \dim Y-2$ (Theorem \ref{thm:scheja2}),
and that in turn gives vanishing results for some classes $H^q(X,\OO(N^{-\mu}))$ (Theorem \ref{thm:scheja3}).
Similarly, we can show easily that
$$H^q(\C\mathbb{P}^k, \OO(N^{-\mu}))=0$$
for all $\mu\geq q-2k$, where $N$ is the universal bundle over $\C\mathbb{P}^k$ (Theorem \ref{thm:extension}).


\section{Behavior of $L^p$-norms under desingularization}\label{sec:resolution}

Let $X$ be a regular irreducible (connected) variety in $\C\mathbb{P}^{n-1}$ of dimension $d-1\geq 1$, \newline
and let $Y$ be the associated homogeneous variety in $\C^n$ (given by the same homogeneous polynomials).
So, $Y$ is an irreducible homogeneous variety in $\C^n$ of dimension $d$, and it is regular outside the origin.
We will now investigate the embedded desingularization of $Y$, which is given by blowing up the origin in $\C^n$.
Let 
$$U \subset \C^n\times \C \mathbb{P}^{n-1}$$
be given by the equations
$$z_j w_k = z_k w_j\ \ \ \mbox{ for all }\ \ j\neq k,$$
where $z_1, ..., z_n$ are the euclidian coordinates of $\C^n$, and $w_1, ..., w_n$ the homogeneous coordinates of $\C\mathbb{P}^{n-1}$.
That is a submanifold of dimension $n$ in $\C^n\times \C \mathbb{P}^{n-1}$.
Let
$$\Pi: U \rightarrow \C^n,\ (z,w) \mapsto z,$$
be the projection to the first component. Then
$$H:= \Pi^{-1}(\{0\}) \cong \C\mathbb{P}^{n-1},$$
but the pre-image of all points in $\C^n\setminus\{0\}$ consists of exactly one point.
We have that
$$\Pi|_{U\setminus H} : U\setminus H \rightarrow \C^n\setminus \{0\}$$
is biholomorphic, $\Pi: U \rightarrow \C^n$ is the blow up of the origin.
On the other hand, consider the projection 
$$P: U \rightarrow H,\ (z,w) \mapsto (0,w).$$
If $\{w_k=1\}$ is a chart in $H$, then 
$$P^{-1}(\{w_k=1\}) \cong \{w_k=1\} \times \C.$$
$U$ is in fact a holomorphic line bundle over $H\cong \C\mathbb{P}^{n-1}$. It is called the universal bundle.
Now, let
$$N := \overline{\Pi_{U\setminus H}^{-1}(Y\setminus\{0\})} \subset U.$$
This is a complex submanifold of dimension $d$ in $U$. Let
$$\pi:=\Pi|_N: N \rightarrow Y,\ \ \ \mbox{ and }\ \  E:=\pi^{-1}(\{0\}) = N\cap H \cong X.$$
Then $\pi: N \rightarrow Y$ is a desingularization of $Y$ (with exceptional set $E\cong X$).
We will from now on identify $E$ with $X$.
On the other hand,
$$p:=P|_N: N \rightarrow X$$
is a holomorphic line bundle. It is the restriction of the universal bundle to $X$, and the normal bundle of $X$ in $N$ at the same time.
Hence, it is a negative bundle in the sense of Grauert (see \cite{Gr1}).
So, there exists an integer $\mu_0\geq 0$ such that
$$H^q(X, \OO(N^{-\mu})) = 0\ \ \ \mbox{ for all }\ \  q\geq 1,\ \mu\geq \mu_0,$$
because the dual bundle $N^{-1} := N^*$ is positive.\\


$U$ is covered by $n$ charts $U_j \cong \C^n$ ($j=1, ..., n$) defined by $w_j=1$.
Let us consider one such domain, say $U_1$. Here, we have holomorphic coordinates
$$z_1,\ w_2,\ ...,\ w_n,$$
and in these coordinates
$$\Pi (z_1, w_2, ..., w_n) = ( z_1, z_1 w_2, ..., z_1 w_n).$$
This implies that
\begin{eqnarray*}
\Pi^* d\o{z_1} &=& d\o{z_1}\ ,\\
\Pi^* d\o{z_j} &=& \o{z_1} d\o{w_j}\ ,\ \ \mbox{ for }\ j=2, ..., n.
\end{eqnarray*}
We will now develop a similar statement on $N$, which is a bit more complicated.
First of all, we will choose a nice hermitian metric $h$ on $U$.
For this, let $h_1'$ be any hermitian metric on $H\cong \C\mathbb{P}^{n-1}$,
say the Fubini-Study metric, and 
$$h_1=P^* h_1'$$
the pull-back to $U$. Furthermore, let $h_2$ be given in the charts $U_j$ (where $w_j=1$) as
$$h_2=\big(|w_1|^2+\cdots +|w_{j-1}|^2 + 1 + |w_{j+1}|^2+\cdots |w_n|^2\big) dz_j\otimes d\o{z_j}.$$
It is easy to see that $h_2$ is globally well defined because $z_j/w_j=z_k/w_k$. Then,
\begin{eqnarray*}
h:=h_1\oplus h_2
\end{eqnarray*}
gives a (in some sense natural) hermitian metric on $U$, where in a chart $U_j$
the coordinate $z_j$ is orthogonal to
$$w_1, ..., w_{j-1}, w_{j+1}, ..., w_n.$$

Let
$$i: N \hookrightarrow U\ \ \ \mbox{ and } \ \ \ \iota: Y \hookrightarrow \C^n$$
be the natural inclusions.
This implies that
$$\Pi\circ i = \iota \circ\pi.$$

As $Y^*$ carries the hermitian structure induced by restriction (respectively pull-back)
of the euclidian metric of the ambient $\C^n$, $N$ is a hermitian submanifold of $U$
with the induced hermitian structure $i^* h$.
We denote by $\|\cdot\|_N$ the resulting norm on the Grassmannian of $N$, and by $dV_N$ 
the associated volume form.\\

Let $Q\in X$ be a point in the exceptional set. We can assume that $Q\in U_1$.
Then there exists a neighborhood $W_Q'$ of $Q$ in $X\cap U_1$ with holomorphic coordinates
$$x_2, ..., x_d$$
on $W_Q'$. It follows that
$$t:=z_1,\ x_2, ..., x_d$$
are holomorphic coordinates on $W_Q:=p^{-1}(W_Q')\subset N$. We identify $x_k$ with $p^* x_k$.
It follows from the construction of the metric that $t=z_1$ is orthogonal to the $x_k$.\\

Hence, by shrinking $W_Q'$ a little, it follows that
\begin{eqnarray}\label{eq:ajk1}
i_* \frac{\partial}{\partial\o{x_j}} = \sum_{k=2}^{n} a_{jk} \frac{\partial}{\partial\o{w_k}},
\end{eqnarray}
where
\begin{eqnarray}\label{eq:ajk2}
\sum_{k=2}^{n} |a_{jk}|^2 \sim 1
\end{eqnarray}
on $W_Q$ for all $j=2, ..., d$. We also have that
\begin{eqnarray*}
dt\wedge d\o{t} \wedge dx_2\wedge d\o{x_2}\wedge \cdots \wedge dx_d\wedge d\o{x_d} \sim dV_N
\end{eqnarray*}
on $W_Q$. 
Using \eqref{eq:ajk1} and \eqref{eq:ajk2}, we calculate:
\begin{eqnarray*}
\big\|\pi_* \frac{\partial}{\partial \o{x_j}}\big\|_Y^2 &=& \big\|\iota_* \pi_* \frac{\partial}{\partial \o{x_j}}\big\|_{\C^n}^2
= \big\| \Pi_* i_*  \frac{\partial}{\partial \o{x_j}}\big\|_{\C^n}^2\\
&=& \big\| \Pi_* \sum_{k=2}^n a_{jk} \frac{\partial}{\partial\o{w_k}}\big\|_{\C^n}^2
= \big\| \sum_{k=2}^n a_{jk} \sum_{l=1}^n \frac{\partial\o{\Pi_l}}{\partial \o{w_k}} \cdot \frac{\partial}{\partial\o{z_l}}\big\|_{\C^n}^2\\
&=& \big\| \sum_{k=2}^n a_{jk} \sum_{l=1}^n \delta_{lk}\cdot z_1 \frac{\partial}{\partial\o{z_l}}\big\|_{\C^n}^2
= |z_1|^2 \sum_{k=2}^n |a_{jk}|^2 \sim |z_1|^2
\end{eqnarray*}
(where $\delta_{lk}$ denotes the Kronecker-$\delta$),
because 
$$\pi_* \frac{\partial}{\partial \o{x_j}}\big|_y\in T^{0,1}_y (Y\setminus\{0\})$$ 
for all $y\in\pi(W_Q)\setminus\{0\}$,
and $\|v\|_Y=\|\iota_* v\|_{\C^n}$ on $T^{0,1} (Y\setminus\{0\})$ (since $\|\cdot\|_Y$ is the norm induced by $\|\cdot\|_{\C^n}$).\\

So, for a point $y\in \pi(W_Q)\setminus\{0\}$, we can now calculate
\begin{eqnarray*}
\big\|(\pi_{N\setminus X}^{-1})^* d\o{x_k}\big\|_Y (y)
&=& \max_{0\neq v\in T^{0,1}_y Y} \|v\|^{-1}_Y(y) \left|d\o{x_k}\big((\pi_{N\setminus X}^{-1})_* v\big)\right|(\pi^{-1}(y))\\
&\sim& \max_{j=2, ..., d} \big\|\pi_* \frac{\partial}{\partial \o{x_j}}\big\|_Y^{-1} (y)
 \left|d\o{x_k}\big(\frac{\partial}{\partial \o{x_j}}\big)\right|(\pi^{-1}(y))\\
 &\sim& \max_{j=2, ..., d} \big\|\pi_* \frac{\partial}{\partial \o{x_j}}\big\|_Y^{-1} (y)
 \sim |z_1(y)|^{-1},
\end{eqnarray*}
because $\pi$ is an biholomorphism outside $X$, and
$$d\o{x_k}\left(\frac{\partial}{\partial\o{t}}\right)=0,$$
since the coordinates $x_2, ..., x_d$ are orthogonal to $t=z_1$.\\

Since $t=\Pi^* z_1 = \pi^* z_1$, the esimate
\begin{eqnarray*}
\big\|(\pi_{N\setminus X}^{-1})^* d\o{x_k}\big\|_Y  \sim |z_1|^{-1}
\end{eqnarray*}
also yields
\begin{eqnarray*}
\big\|(\pi_{N\setminus X}^{-1})^* ( \o{t} d\o{x_j})\big\|_Y  \sim 1.
\end{eqnarray*}
Summing up, we conclude:

\begin{lem}\label{lem:blowup}
Let $Q\in X$. Then there exists a neighborhood $W_Q$ of $Q$ in $N$
with holomorphic coordinates $t, x_2, ..., x_d$ such that
$$X\cap W_Q = \{ t=0\},$$
and
\begin{eqnarray*}
\alpha_1 &:=& (\pi|_{N\setminus X}^{-1})^* d\o{t},\\
\alpha_j &:=& (\pi|_{N\setminus X}^{-1})^* (\o{t} d\o{x_j}),\ \ j=2, ..., d,
\end{eqnarray*}
are a basis of the $(0,1)$-forms on $\pi(W_Q \setminus X)\subset Y\setminus\{0\}$ with
$$\|\alpha_j\|_Y \sim 1.$$
This implies for the volume forms that
$$\pi|_{N\setminus X}^* dV_Y \sim |t|^{2d-2} dV_N$$
on $W_Q\setminus X$. Hence, for a function $f$ on $\pi(W_Q \setminus X)=\pi(W_Q)\setminus\{0\}$, and $1\leq p \leq \infty$,
we have that
$$f\in L^p(\pi(W_Q) \setminus \{0\})$$
exactly if
$$|t|^{\frac{2d-2}{p}} \cdot \pi^* f \in L^p(W_Q \setminus X).$$
Let $1\leq q\leq n$. If $\omega \in L^p_{0,q}(\pi(W_Q)\setminus \{0\})$ is a $(0,q)$-form on $\pi(W_Q) \setminus\{0\}$, then
\begin{eqnarray*}
|t|^{\frac{2d-2}{p} - (q-1)} \cdot \pi^* \omega \in L^p_{0,q}(W_Q\setminus X).
\end{eqnarray*}
On the other hand, for $\eta\in L^p_{0,q}(W_Q\setminus X)$ a $(0,q)$-form on $W_Q\setminus X$,
$$|t|^{\frac{2d-2}{p} - q} \cdot \eta \in L^p_{0,q}(W_Q\setminus X)$$
implies that
\begin{eqnarray*}
(\pi|_{N\setminus X}^{-1})^* \eta \in L^p_{0,q}(\pi(W_Q)\setminus \{0\}).
\end{eqnarray*}
\end{lem}

\begin{proof}
Only the last two statements remain to show. $\omega\in L^p_{0,q}(\pi(W_Q)\setminus\{0\})$ has a representation
$$\sum_{1\leq k_1 < \cdots < k_q\leq d} f_{k_1\cdots k_q} \alpha_{k_1}\wedge\cdots \wedge \alpha_{k_q},$$
where the coefficients $f_{k_1\cdots k_d} \in L^p(\pi(W_Q)\setminus\{0\})$, and the proof is clear from what we have seen before.
The last statement follows analogously.
\end{proof}


\section{Sufficient Conditions (Theorem \ref{thm:sufficient})}\label{sec:sufficient}

Let $D$ be a strongly pseudoconvex domain in $Y$ such that $0\in D$,
and let
$D^*:=D\setminus\{0\}$.
We can assume that $D\cap U =\{z\in U:\rho(z)<0\}$ where $\rho\in C^2(U)$ is a regular strictly plurisubharmonic 
defining function on a neighborhood $U$ of $bD$.
Then there exists $\epsilon>0$ such that $D_\epsilon:= D\cup \{z\in U:\rho(z)<\epsilon\}$ is a strongly pseudoconvex extension of $D$.
So, it follows by Grauert's bump method that the natural homomorphism
$$r_q: H^q_{(p)}(D_\epsilon^*,\OO) \rightarrow H^q_{(p)}(D^*,\OO)$$
(induced by restriction of forms) is surjective (see \cite{LiMi}, chapter IV.7).
Here, we also set $D_\epsilon^*=D_\epsilon\setminus\{0\}$.
We will work with the desingularization 
$\pi: N \rightarrow Y$
described in the previous section. So, let
$$G=\pi^{-1}(D),\ G_\epsilon=\pi^{-1}(D_\epsilon),\ G^*=\pi^{-1}(D^*),\ G^*_\epsilon=\pi^{-1}(D^*_\epsilon),$$
and $[\omega] \in H^q_{(p)}(D^*,\OO)$ represented by $\omega\in L^p_{0,q}(D^*_\epsilon)$.
We will show in this section how $\omega$ determines a class in \eqref{eq:iso},
and that $[\omega]=0$ if that class vanishes.
The point that a different representative of $[\omega]$ defines the same class
is postponed to the next section.
We can use Lemma \ref{lem:blowup} to determine properties of $\pi^* \omega$.
It is convenient to work with the weighted Dolbeault complexes that we introduced in \cite{Rp6}.
So, we have to describe some concepts.
Let $\mathcal{I}$ be the sheaf of ideals of $E=X$ in $N$. For $k\in\Z$ we will use
the sheaves $\mathcal{I}^k\OO$ which are subsheaves of the sheaf of germs of meromorphic functions on $N$.
It follows from Theorem 5.1 in \cite{Rp4} that
\begin{eqnarray}\label{eq:iso}
H^q(G,\mathcal{I}^k\OO) \cong H^q(G_\epsilon,\mathcal{I}^k\OO) \cong \bigoplus_{\mu\geq k} H^q(X,\OO(N^{-\mu}))
\end{eqnarray}
for all $q\geq 1$, because $G$ and $G_\epsilon$ are strongly pseudoconvex neighborhoods of the zero section of the
negative holomorphic line bundle $N$.
Note that on the left-hand side of \eqref{eq:iso}, $\OO$ is the structure sheaf on $N$,
while on the right $\OO(N^{-\mu})$ denotes the sheaf of germs of holomorphic sections in the bundle $N^{-\mu}$
over $X$.
So, in order to prove Theorem \ref{thm:sufficient}, it is enough to show that there exists an injective
homomorphism
\begin{eqnarray}\label{eq:sufficient2}
H^q_{(p)}(D^*,\OO) \hookrightarrow H^q(G,\mathcal{I}^{a(p,q,d)}\OO).
\end{eqnarray}
What we need is a suitable fine resolution for the sheaves $\mathcal{I}^k\OO$.
Let $s\in \R$, $U\subset N$ open,
and $\eta$ a measurable $(0,r)$-form on $U$.
Then, we say that
$$\eta\in |\mathcal{I}|^s L^p_{(0,r),loc}(U)$$
if for each point $z\in U$ there is a local generator $f_z$ of $\mathcal{I}_z$ (defined on a neighborhood $V_z$ of $z$)
such that
$$|f_z|^{-s} \eta \in L^p_{0,r}(V_z).$$
This property does not depend on the choice of $f_z$, and so the spaces $|\mathcal{I}|^s L^p_{(0,r),loc}(U)$
are well-defined.\\

We have to use a weighted $\dq$-operator, which we define locally again.
Let $k\in\Z$,
$z\in N$ and $f_z$ a local generator of $\mathcal{I}_z$ defined on $V_z$.
Then, for a current $\Phi$ on $V_z$, we set
$$\dq_k \Phi:= f_z^k \dq \big( f_z^{-k} \Phi\big),$$
provided the construction makes sense.
In that case $\dq_k$ is well-defined because the construction does not depend on the choice of the generator.
Now, we have to make a connection between the weighted operators $\dq_k$
and weighted $L^p$-spaces defined above. We will use:

\begin{defn}\label{defn:k}
Let $1\leq p \leq \infty$ and $s$ be real numbers. Then we call
\begin{eqnarray*}
k(p,s) := \max\{m\in\Z: |z_1|^s L^p_{loc}(\C) \subset |z_1|^m L^1_{loc}(\C)\}
\end{eqnarray*}
the $\dq$-weight of $(p,s)$, where $|z_1|^t L^p_{loc}(\C)=\{f \mbox{ measurable}: |z_1|^{-t} f \in L^p_{loc}(\C)\}$.
Now, we define for $0\leq q\leq d=\dim Y$ the sheaves $|\mathcal{I}|^s \mathcal{L}^p_{0,r}$ by:
\begin{eqnarray*}
|\mathcal{I}|^s \mathcal{L}^p_{0,r} (U) :=\{ f\in |\mathcal{I}|^s L^p_{(0,r),loc}(U): \dq_{k(p,s)} f \in |\mathcal{I}|^s L^p_{(0,r+1),loc}(U)\}
\end{eqnarray*}
for open sets $U\subset\C^n$ (it is a presheaf wich is already a sheaf).
\end{defn}

From now on, if an index $k$ is not specified, it should always
be the $\dq$-weight $k(p,s)$, where $p$ and $s$ arise from the context.
We need to compute the $\dq$-weight of $(p,s)$ explicitly:

\begin{lem}\label{lem:k1}
Let $1\leq p \leq \infty$ and $s$ be real numbers, and $k(p,s)$ the $\dq$-weight of $(p,s)$
according to Definition \ref{defn:k}. Then
\begin{eqnarray}\label{eq:k1}
k(p,s) = \left\{
\begin{array}{ll}
\max\{m\in\Z: m<2 + s-2/p\} & ,\ p\neq 1,\\
\max\{m\in\Z: m\leq 2 + s-2/p\}& ,\ p=1.
\end{array}\right.
\end{eqnarray}
\end{lem}

\begin{proof}
See \cite{Rp6}, Lemma 2.2.
\end{proof}

We can now cite the main results about the Dolbeault complex with weights according to normal crossings.
Adapted to our present situation, Theorem 1.5 in \cite{Rp6} reads as:

\begin{thm}\label{thm:main}
For $1\leq p\leq \infty$ and $s\in \R$, let $k(p,s)\in \Z$ be the $\dq$-weight according to Definition \ref{defn:k}.
Then:
\begin{eqnarray}\label{eq:complex}
0 \rightarrow \mathcal{I}^k\OO \hookrightarrow
 |\mathcal{I}|^s \mathcal{L}^p_{0,0} \xrightarrow{\ \dq_k\ }
 |\mathcal{I}|^s \mathcal{L}^p_{0,1} \xrightarrow{\ \dq_k\ }
\cdots \xrightarrow{\ \dq_k\ }
 |\mathcal{I}|^s \mathcal{L}^p_{0,d} \rightarrow 0
\end{eqnarray}
is an exact (and fine) resolution of $\mathcal{I}^k\OO$.
\end{thm}

Let us now return to $\omega\in L^p_{0,q}(D^*_\epsilon)$.
If we extend $\pi^*\omega$ trivially over the exceptional set $E$,
Lemma \ref{lem:blowup} implies immediately:

\begin{lem}\label{lem:piw}
Let $s=(q-1) - \frac{2d-2}{p}$. Then:
$$\pi^* \omega \in |\mathcal{I}|^s L^p_{(0,q),loc} (G_\epsilon).$$
\end{lem}

Now, we need to find a suitable weight $k$ such that $\dq_k \pi^*\omega=0$.
This is in fact the $\dq$-weight of $(p,s)$ in Lemma \ref{lem:piw},
as we will see shortly. But before, it is the time to make the connection to Theorem \ref{thm:sufficient}:

\begin{lem}\label{lem:piw1}
Let $k(p,s)$ be the $\dq$-weight of $p$ and 
$$s=(q-1) - \frac{2d-2}{p}.$$
Then:
$$k(p,s) = a(p,q,d),$$
where $a(p,q,d)$ is the constant from Theorem \ref{thm:sufficient}.
So, Lemma \ref{lem:piw} yields
\begin{eqnarray}\label{eq:al1}
\pi^* \omega \in |\mathcal{I}|^s L^p_{(0,q),loc}(G_\epsilon) \subset |\mathcal{I}|^{a(p,q,s)} L^1_{(0,q),loc} (G_\epsilon)
\end{eqnarray}
by Definition of the $\dq$-weight.
\end{lem}

\begin{proof}
The proof is immediate, because
\begin{eqnarray*}
2+s-2/p &=& (q+1) - \frac{2d-2}{p} - 2/p = (q+1) -2d/p.
\end{eqnarray*}
\end{proof}

From now on, if the indices are not specified, $a$ should always be the constant $a(p,q,d)$ from Theorem \ref{thm:sufficient}.
We will now see that in fact
$$\dq_a \pi^*\omega=0.$$ 
This is a consequence of \eqref{eq:al1}, Lemma \ref{lem:blowup} and
the following extension theorem for the $\dq$-equation, which we will show in a (for further use) slightly more general version than needed:

\begin{lem}\label{lem:ext1}
Let $D\subset\C^n$ be an open set, $1\leq P\leq\infty$ and $f\in L^P_{0,Q}(D)$
a $(0,Q)$-form on $D$ such that $\dq f =g$ in the sense of distributions
on $D\setminus H$, where $H=\{z\in \C^n: z_1=0\}$ and $g\in L^1_{0,Q+1}(D)$, and $f$ has the following structure:
$$ f =\sum_{|J|=Q} f_J d\o{z_J}$$
(in multi-index notation) such that 
$$|z_1|^{-w(P)} f_J \in L^P(D) \mbox{ for all multi-indices } J \mbox{ with } 1\notin J,$$
where
$$w(P)=\left\{\begin{array}{ll}
2/P -1 & , \mbox{ if } 1\leq P \leq 2,\\
0& , \mbox { if } 2\leq P \leq \infty.\end{array}\right.$$
Then $\dq f=g$ on the whole set $D$.
\end{lem}

We will use the statement only in case $P=1$ and $w(P)=1$.

\begin{proof}
The statement is local, so we can assume that $D$ is bounded.
For $r>0$, define
$$U(r):=\{z\in\C^n: \dist(z,H)=|z_1|<r\}.$$
Choose a smooth cut-off function $\chi\in C^\infty_{cpt}(\R)$ with $|\chi|\leq 1$,
$\chi(t)=1$ if $|t|\leq 1/2$, $\chi(t)=0$ if $|t|\geq 2/3$, and $|\chi'| \leq 8$.
Now, let
$$\chi_r(z):=\chi(\frac{\dist(z,H)}{r}).$$
Then $\chi_r\equiv 1$ on $U(r/2)$ and $\supp \chi_r \subset U(3r/4)$.
$\chi_r$ is smooth
and we have:
$$\||z_1|^{w}\nabla \chi_r\|\leq |\chi'| r^{w-1} \leq 8r^{w-1}.$$
Since $D$ is bounded, there is $R>0$ such that $D\subset B_R(0)$. 
Let $s=P/(P-1)$ be the coefficient dual to $P$. 
It follows that
\begin{eqnarray*}
\int_{D}\||z_1|^w \nabla \chi_r\|^s dV_{\C^n} 
\leq 8^s (2R)^{2n-2} r^{s(w-1)} \int_{\{\zeta\in\C:|\zeta|<r\}} dV_{\C},
\end{eqnarray*}
and we conclude:
\begin{eqnarray}\label{eq:dchi}
\||z_1|^w \nabla\chi_r\|_{L^s(D)} &\lesssim& r^{w-1 + 2/s} = r^{w-1+2 - 2/P}\lesssim 1
\end{eqnarray}
by the choice of $w(P)$. The statement remains true in case $P=1$ and $s=\infty$.
What we have to show is that
\begin{eqnarray}\label{eq:dqe2}
\int_D f\wedge \dq \phi = (-1)^{q+1} \int_D g\wedge\phi
\end{eqnarray}
for all smooth $(n,n-Q-1)$-forms $\phi$ with compact support in $D$.
By assumption, $\dq f=g$ on $D\setminus H$. That leads to:
\begin{eqnarray*}
\int_D f\wedge\dq\phi &=& \int_D f\wedge\chi_r \dq\phi +\int_{D\setminus H} f \wedge(1-\chi_r) \dq\phi\\
&=& \int_D f \wedge\chi_r \dq\phi + \int_{D\setminus H} f \wedge\dq[ (1-\chi_r) \phi]-
\int_{D\setminus H} f\wedge\dq(1-\chi_r)\wedge \phi\\
&=& \int_D f \wedge\chi_r \dq\phi + (-1)^{Q+1} \int_{D} g \wedge(1-\chi_r) \phi 
+ \int_{D\setminus H} f \wedge\dq\chi_r\wedge\phi.
\end{eqnarray*}
Now, we will consider what happens as $r\rightarrow 0$.
Let us first consider
$$\int_D f\wedge\chi_r\dq\phi\ \ \mbox{ and }\ \ \int_{D} g \wedge(1-\chi_r) \phi.$$
Since $|\chi_r|\leq 1$, we have
\begin{eqnarray*}
\|f\wedge\chi_r\dq\phi \|, \|g\wedge(1-\chi_r)\phi\|  \in L^1(D),
\end{eqnarray*}
and we know that $f\wedge\chi_r\dq\phi \rightarrow 0$ pointwise if $r\rightarrow 0$,
and $g\wedge (1-\chi_r)\phi \rightarrow g\wedge\phi$.
Hence, Lebesgue's Theorem on dominated convergence gives:
\begin{eqnarray*}
\lim_{r\rightarrow 0} \int_D f \wedge\chi_r \dq\phi = 0,\ \ \ \lim_{r\rightarrow 0} \int_D g \wedge(1- \chi_r) \phi = \int_D g\wedge\phi.
\end{eqnarray*}

To prove \eqref{eq:dqe2}, only
$$\lim_{r\rightarrow 0}\int_{D} f\wedge\dq\chi_r\wedge\phi=0$$
remains to show. Because of
$$\dq \chi_r = \frac{\partial\chi_r}{\partial\o{z_1}} d\o{z_1},$$
we only have to consider the coefficients $f_J$ where $1\notin J$. So,
using \eqref{eq:dchi} and the H\"older Inequality, we get
\begin{eqnarray*}
\lim_{r\rightarrow 0} \|f\wedge\dq\chi_r\wedge \phi\|_{L^1(D)}
&=& \lim_{r\rightarrow 0} \|f\wedge\dq\chi_r\wedge \phi\|_{L^1(U(r))}\\
&\leq& \lim_{r\rightarrow 0} \|f \wedge\phi\|_{L^P(U(r))}\||z_1|^w \nabla \chi_r\|_{L^s(D)}\\
&\lesssim&  \lim_{r\rightarrow 0} \|f \wedge\phi\|_{L^P(U(r))}.
\end{eqnarray*}
Since $f\in L^P$, we conclude
$$\lim_{r\rightarrow 0} \|f \wedge\phi\|_{L^P(U(r))} = 0$$
(see for instance \cite{Alt}, Lemma A 1.16), and that completes the proof.
\end{proof}

So, choose a point on the exceptional set $E$.
Locally, we can assume that this point is the origin in $\C^d$,
and that $E=\{z_1=0\}$ in a small neighborhood $V$. It follows from Lemma \ref{lem:piw1} that
$$z_1^{-a} \pi^* \omega \in L^1_{0,q}(V),$$
and it is clear that 
$$\dq \big( z_1^{-a} \pi^* \omega\big) =0 \ \mbox{ on }\ V\setminus E$$
in the sense of distributions. But if we take a closer look at $z_1^{-a} \pi^*\omega$
it follows from Lemma \ref{lem:blowup}, that
$$z_1^{-a} \pi^*\omega = \sum_{|J|=q} f_J d\o{z_J},$$
where $f_J=|z_1| h_J$ with $h_J\in L^1(V)$ if $1\notin J$.
So, Lemma \ref{lem:ext1} yields $\dq (z_1^{-a} \pi^*\omega)=0$ on $V$.
Thinking globally, that means nothing else but:

\begin{lem}\label{lem:piw2}
If $a=a(p,q,d)$ is the index from Theorem \ref{thm:sufficient}, then
$$\dq_a \pi^* \omega=0.$$
Hence, $\pi^* \omega\in |\mathcal{I}|^s \mathcal{L}^p_{0,q}(G_\epsilon)$ (where $s=(q-1)-\frac{2d-2}{p}$) defines a cohomology class
$$[\pi^*\omega] \in H^q(G_\epsilon, \mathcal{I}^a \OO).$$
\end{lem}

\begin{proof}
Because of $a(p,q,d)=k(p,s)$ by Lemma \ref{lem:piw1},
$\dq_a \pi^*\omega=\dq_k \pi^*\omega=0$ implies that $\pi^* \omega\in |\mathcal{I}|^s \mathcal{L}^p_{0,q}(G_\epsilon)$.
Hence, $\pi^*\omega$ defines a cohomology class $[\pi^*\omega]$ in $H^q(G_\epsilon,\mathcal{I}^a\OO)$ by Theorem \ref{thm:main}.
\end{proof}

Now, assume that
\begin{eqnarray}\label{eq:closed1}
[\pi^*\omega]=0\ \ \ \mbox{ in } \ \ \ H^q(G_\epsilon,\mathcal{I}^a\OO).
\end{eqnarray}
We will conclude this section by showing that this implies $[\omega]=0$ in $H^q_{(p)}(D^*,\OO)$.
By the use of Theorem \ref{thm:main}, the assumption \eqref{eq:closed1} tells us that there exists
\begin{eqnarray}\label{eq:closed2}
\eta\in |\mathcal{I}|^s L^p_{(0,q-1),loc}(G_\epsilon)
\end{eqnarray}
such that $\dq_k\eta=\pi^*\omega$ on $G_\epsilon$. This means that $\dq\eta=\pi^*\omega$ on $G_\epsilon \setminus E=G_\epsilon^*$.
Recall that
\begin{eqnarray}\label{eq:closed3}
s=(1-q) -\frac{2d-2}{p}.
\end{eqnarray}
Let $\eta':=\eta|_{G}$.
Then, \eqref{eq:closed2}, \eqref{eq:closed3} and the last statement of Lemma \ref{lem:blowup}
yield that
$$\theta := (\pi|_{G \setminus E}^{-1})^* \eta \in L^p_{0,q-1}(D\setminus \{0\}).$$
Because $\pi$ is a biholomorphic map outside the exceptional set,
we know that $\dq\theta=\omega$ on $D\setminus \{0\}=D^*$.
So, it follows that $[\omega]=0$ in $H^q_{(p)}(D^*,\OO)$.
To complete the proof of Theorem \ref{thm:sufficient},
it remains to show that a different representing $(0,q)$-form
for the class
$[\omega]\in H^q_{(p)}(D^*,\OO)$
defines the same class in
$$H^q(G_\epsilon,\mathcal{I}^{a}\OO) \cong H^q(G,\mathcal{I}^a\OO).$$
That will be done in the next section, where we can restrict our considerations to the set $G$
(no need to consider the extension $G_\epsilon$ any more).

\vspace{2mm}
\section{Integration along the Fibers}\label{sec:fibers}

Assume that $\wt{\omega}$ is another representing form for the class $[\omega]\in H^q_{(p)}(D^*,\OO)$,
namely that 
$$\wt{\omega} \in L^p_{0,q}(D^*),$$
such that $\dq \wt{\omega}=0$ on $D^*$ in the sense of distributions, and
there exists 
$\sigma \in L^p_{0,q-1}(D^*)$
with
$$\omega-\wt{\omega} = \dq \sigma.$$
Here again, $\pi^*\wt{\omega}\in |\mathcal{I}|^s \mathcal{L}^p_{0,q}(G)$,
but unfortunately we do not have $\pi^* \sigma \in |\mathcal{I}|^s\mathcal{L}^p_{0,q-1}(G)$.
But, we can use $\pi^*\sigma$ to construct $\wt{\sigma}\in |\mathcal{I}|^s\mathcal{L}^p_{0,q-1}(G)$
such that 
$$\pi^*\omega-\pi^*\wt{\omega} = \dq_a \wt{\sigma}.$$
Let $\chi \in C^\infty_{cpt}(G)$ be a smooth cut-off function with compact support in $G$ such that $\chi\equiv 1$
in a neighborhood of the exceptional set $E=X$. Now, consider
\begin{eqnarray*}
\delta := \pi^* \omega-\pi^*\wt{\omega} -\dq \big((1-\chi)\pi^*\sigma\big) \in |\mathcal{I}|^s \mathcal{L}^p_{0,q}(G).
\end{eqnarray*}
We will now solve the equation $\dq_a\tau=\delta$ with $\tau\in |\mathcal{I}|^s\mathcal{L}^p_{0,q-1}(G)$,
and then
$$\dq_a\big( \tau + (1-\chi)\pi^* \sigma \big) = \pi^*\omega - \pi^*\wt{\omega}$$
will tell us that $[\pi^*\omega] = [\pi^*\wt{\omega}]$ in $H^q(G,\mathcal{I}^a \OO)$.

The crucial point is that the form $\delta$ has compact support in $G$.
That allows us to solve the equation $\dq_a\tau=\delta$ by integrating over the fibers of $N$ interpreted as a holomorphic line bundle
over $E=X$. That idea has been already used by E.\ S.\ Zeron and the author in \cite{RuZe}.
We can define that integration locally: Let $Q\in X$.
Then there exists a neighborhood $U_Q$ of $Q$ in $X$ with coordinates $z_1, ..., z_{d-1}$ such that
$N$ is trivial over $U_Q$:
$$N|_{U_Q} \cong U_Q\times \C.$$
So, let $z_1, ... z_{d-1}, t$ be the coordinates on $V_Q:=N|_{U_Q}$.
Then, $\delta|_{V_Q} \in |t|^s L^p_{0,q}(V_Q)$ can be written uniquely as
\begin{eqnarray*}
\delta|_{V_Q} = \sum_{|J|=q} g_J d\o{z_J} + \sum_{|J|=q-1} f_{J} d\o{t}\wedge d\o{z_J}
\end{eqnarray*}
where all the coefficients $g_J, f_{J} \in |t|^s L^p(V_Q)$ satisfy $t^{-a} g_J, t^{-a} f_J \in L^1(V_Q)$.
Now, we define:
\begin{eqnarray}\label{eq:fibers}
\tau|_{V_Q} := \sum_{|J|=q-1} \left( \frac{t^a}{2\pi i} \int_{\C} \frac{f_J(z_1, ..., z_{d-1}, \zeta)}{\zeta^a} \frac{d\zeta \wedge d\o{\zeta}}{\zeta-t} \right) d\o{z_J}.
\end{eqnarray}
We have to show that this construction globally defines a form $\tau\in |\mathcal{I}|^s \mathcal{L}^p_{0,q-1}(G)$ such that $\dq_a\tau=\delta$,
where we intensively use the fact that $\delta$ has compact support.
Firstly, we remark that the operator in \eqref{eq:fibers} maps continuously
$$|t|^s L^p_{0,q}(V_Q\cap G) \rightarrow |t|^{s+1-\epsilon} L^p_{0,q-1}(V_Q\cap G) \subset |t|^s L^p_{0,q-1}(V_Q\cap G),$$
because $a(p,q,d)=k(p,s)$ is the $\dq$-weight of $(p,s)$ (see \cite{Rp6}, Theorem 2.1).
For $\dq_a \tau=\delta$, we have to show that 
\begin{eqnarray}\label{eq:dqt}
\dq\big( t^{-a} \tau|_{V_Q}\big) = t^{-a} \delta|_{V_Q}.
\end{eqnarray}
Here, we have to work a little, because we are dealing with weak concepts.
Let $U\subset \C^{d-1}$ be an open set, 
$1\leq r \leq d$, and $\omega\in C^0_{0,r}(U\times\C)$
such that 
$\omega$ has compact support in the $z_d$-direction,
i.e. $\supp \omega\cap F_a$ is compact in $F_a$ for all fibers
$F_a=\{a\}\times\C$, $a\in U$. If
$$\omega = \sum_{\substack{|J|=r-1,\\ d\notin J}} a_{dJ} d\o{z_d}\wedge d\o{z_J} 
+ \sum_{\substack{|K|=r,\\ d\notin K}} a_K d\o{z_K}$$
(the multi-indices in ascending order), then let
$${\bf S}_r^d \omega := \sum_{|J|=r-1} {\bf I}(a_{dJ})\ d\o{z_J}$$
where
$${\bf I} f (z_1, ..., z_d) :=\frac{1}{2\pi i} \int_\C f(z_1, ..., z_{d-1}, t) \frac{dt \wedge d\o{t}}{t-z_d}.$$


It is not hard to see that we can use these operators ${\bf S}_r^d$ to construct a $\dq$-homotopy formula
for forms with compact support in the fibers:

\begin{lem}\label{lem:dqc}
Let $\omega\in C^1_{0,q}(U\times\C)$ such that $\omega$ has compact support in $z_d$.
Then:
$$\omega = \dq {\bf S}^d_q \omega + {\bf S}^d_{q+1}\dq \omega.$$
\end{lem}

\begin{proof}
Let
$$\omega = \sum_{\substack{|J|=q-1,\\ d\notin J}} a_{dJ} d\o{z_d}\wedge d\o{z_J} 
+ \sum_{\substack{|K|=q,\\ d\notin K}} a_K d\o{z_K}$$
and
$$\dq \omega = \sum_{\substack{|K|=q,\\ d\notin K}} c_{dK} d\o{z_d}\wedge d\o{z_K} + \cdots$$
Then we compute that:
\begin{eqnarray}\label{eq:dq-closed}
c_{dK}=\frac{\partial a_K}{\partial \o{z_d}} - \sum_{\substack{J\subset K,\\ |J|=q-1}} \mbox{sign} { K\setminus J\ J \choose K} 
\frac{\partial a_{dJ}}{\partial \o{z_{K\setminus J}}}.
\end{eqnarray}
By use of the inhomogeneous Cauchy-Integral Formula
in one complex variable and the assumption about the support of $\omega$,
we compute furthermore:
\begin{eqnarray*}
\dq \big({\bf I} (a_{dJ}) d\o{z_J}\big)
= a_{dJ} d\o{z_d}\wedge d\o{z_J} + \sum_{k\notin J\cup \{d\}} {\bf I}\left(\frac{\partial a_{dJ}}{\partial \o{z_k}}\right)
d\o{z_k}\wedge d\o{z_J}.
\end{eqnarray*}
But this leads to (summing up):
$$\dq {\bf S}_q^d \omega 
=
\sum_{\substack{|J|=q-1,\\ d\notin J}} a_{dJ} d\o{z_d}\wedge d\o{z_J} 
+ \sum_{\substack{|K|=q,\\ d\notin K}} {\bf I}\big(b_K\big) d\o{z_K},$$
where
\begin{eqnarray*}
b_K &=& 
\sum_{J\subset K} \mbox{sign} { K\setminus J\ J \choose K} 
\frac{\partial a_{dJ}}{\partial \o{z_{K\setminus J}}} = \frac{\partial a_K}{\partial \o{z_d}} - c_{dK}
\end{eqnarray*}
by the use of \eqref{eq:dq-closed}. But $a_K$ has compact support in $z_d$.
So (using the inhomogeneous Cauchy Formula again),
\begin{eqnarray*}
{\bf{I}}\left(b_K\right) d\o{z_K} &=& {\bf I}\left(\frac{\partial a_K}{\partial \o{z_d}}\right) d\o{z_K} - {\bf I}\big( c_{dK} \big) d\o{z_K}\\
&=& a_K d\o{z_K} - {\bf I}\big( c_{dK} \big) d\o{z_K} ,
\end{eqnarray*}
and we are done, because
$${\bf S}_{q+1}^d \dq\omega =\sum_{\substack{|K|=q,\\d\notin K}} {\bf I}\big( c_{dK}\big) d\o{z_K}.$$
\end{proof}

Turning to $L^1$-forms, we can deduce:

\begin{lem}\label{lem:dqc2}
Let $R>0$ and $\omega\in L^1_{0,q}(U\times\C)$ such that $\dq\omega=0$, and $\omega$ has support in $U\times \Delta_R(0)$,
where $\Delta_R(0)$ is the disc of radius $R$ at $0$.
Then:
$$\omega = \dq {\bf S}^d_q \omega$$
on each subset $V\times \Delta_R(0)$ where $V\subset\subset U$.
\end{lem}

\begin{proof}
We simply use the assumption $V\subset\subset U$ because we really do not need to care about the boundary.
Using convolution with a Dirac sequence, there exists a sequence of smooth forms $f_j \in C^\infty_{0,q}(V\times \C)$
such that 
\begin{eqnarray*}
\lim_{j\rightarrow \infty} f_j &=& \omega \ \ \mbox{ in }\ L^1_{0,q}(V\times \C),\\
\lim_{j\rightarrow \infty} \dq f_j &=& 0 \ \ \mbox{ in }\ L^1_{0,q+1}(V\times \C),
\end{eqnarray*}
and we can assume that the $f_j$ have support in $V\times \Delta_{R+1}(0)$. Lemma \ref{lem:dqc} tells us that
$$f_j = \dq {\bf S}^d_q f_j + {\bf S}^d_{q+1} \dq f_j$$
for all $j$, and passing to the limit in $L^1$-spaces proves the Lemma,
because the operators ${\bf S}^d_q$, ${\bf S}^d_{q+1}$ map continuously from $L^1$ to $L^1$.
\end{proof}

Let us return to the $\dq$-equation \eqref{eq:dqt} which we are trying to prove.
But that is now an easy consequence of Lemma \ref{lem:dqc2},
because
$$t^{-a} \tau|_{V_Q} = {\bf S}^d_q \big(t^{-a} \delta|_{V_Q}\big).$$

It only remains to show that $\tau$ is globally well-defined.
If we change coordinates on $X$, that does not effect the Definition \eqref{eq:fibers},
but we have to care about what happens for a different trivialization of $N$.
So, let $w=\phi(z_1, ..., z_{d-1}) t$ and 
\begin{eqnarray*}
\delta|_{V_Q} = \sum_{|J|=q} g_J d\o{z_J} + \sum_{|J|=q-1} \wt{f_{J}} d\o{w}\wedge d\o{z_J}.
\end{eqnarray*}
Then $d\o{w} = \o{\phi(z_1, ..., z_{d-1})} d\o{t}$ and 
$$\wt{f_J} = \o{\phi}^{-1} f_J.$$
That yields (with $\xi=\phi\zeta$):
\begin{eqnarray*}
\frac{t^a}{2\pi i} \int_{\C} \frac{f_J(z_1, ..., z_{d-1},\zeta)}{\zeta^a} \frac{d\zeta \wedge d\o{\zeta}}{\zeta-t} &=&
\frac{\phi^{-a} w^a}{2\pi i} \int_{\C} \frac{\o{\phi} \wt{f_J}(z_1, ..., z_{d-1},\xi)}{\phi^{-a}\xi^a} \frac{|\phi|^{-2} d\xi \wedge d\o{\xi}}{\phi^{-1}(\xi-w)}\\
&=& \frac{w^a}{2\pi i} \int_{\C} \frac{\wt{f_J}(z_1, ..., z_{d-1},\xi)}{\xi^a} \frac{d\xi \wedge d\o{\xi}}{\xi-w},
\end{eqnarray*}
and that shows that $\tau \in |\mathcal{I}|^s L^p_{0,q-1}(G)$ is globally well defined by \eqref{eq:fibers}.
Since $\dq_a \tau=\delta$ as we have seen before, we have finished the proof that $[\pi^* \omega]=[\pi^*\wt{\omega}]$ in $H^q(G,\mathcal{I}^a \OO)$,
and that also finishes the proof of Theorem \ref{thm:sufficient}.


Another interesting application of the integration along the fibers is Theorem \ref{thm:integration},
namely the solution of the $\dq$-equation for forms with compact support:

\vspace{3mm}
{\bf Theorem 1.3.}
{\it Let $X$ and  $Y$ be as in the introduction, $D\subset\subset Y$ an open subset, $D^*=D\setminus\{0\}$,
and $1\leq p\leq \infty$, $1\leq q\leq \dim Y$. Let
$$\omega \in L^p_{0,q}(D^*)\cap \ker \dq$$
with compact support in $D$. Then there exists
$\eta\in L^p_{0,q-1}(D^*)$
such that $\dq\eta=\omega$.}

\begin{proof}
As before, let $G:=\pi^{-1}(D)$ and $G^*:=G\setminus E$.
Let $\omega\in L^p_{0,q}(D^*)$ such that $\dq\omega=0$, and $\omega$ has support in $D$.
Then (by Lemma \ref{lem:ext1})
$$\pi^*\omega \in |\mathcal{I}|^s L^p_{0,q}(G),$$
where 
$$s=(q-1)-\frac{2d-2}{p},$$
and by Lemma \ref{lem:ext1} we have $\dq_a\pi^*\omega=0$ with $a=a(p,q,d)=k(p,s)$
the $\dq$-weight of $(p,s)$.
Because $\pi^* \omega$ has compact support in $G$, we can integrate along the fibers as before
and obtain $\tau \in |\mathcal{I}|^s L^p_{0,q-1} (G)$
such that $\dq_a \tau=\pi^*\omega$. But then
$$\eta:= (\pi|_{G\setminus E}^{-1})^* \tau \in L^p_{0,q-1}(D^*)$$
by Lemma \ref{lem:blowup}, and 
$\dq \eta= \omega$ on $D^*$.
\end{proof}

As a consequence, we can now also prove Theorem \ref{thm:compact}, namely $H^1_{(p),cpt}(D^*,\OO)=0$
for $D\subset\subset Y$ as above. So, let $\omega\in L^p_{0,1}(D^*)$ be $\dq$-closed with support in $D$.
We can assume that
$$D\subset\subset \wt{D}=Y\cap B_R(0)$$
for $R>0$ large enough, and extend $\omega$ trivially by $0$ to the whole set $\wt{D}^*=\wt{D}\setminus\{0\}$.
By Theorem \ref{thm:integration}, there exists $f\in L^p(\wt{D}^*)$ such that $\dq f= \omega$.
So, $f$ is holomorphic on $\wt{D}\setminus D$ which is connected because $X$ and $Y$ are chosen irreducible.

\vspace{2mm}
But $Y$ is a normal complex space, because $Y$ is a complete intersection, and
a complete intersection is a normal space exactly if the codimension of the singular set is $\geq 2$
(see \cite{Ab}, 12.3, or \cite{Sch2}, Korollar 4).

\vspace{2mm}
So, $f|_{\wt{D}\setminus D}$ extends uniquely to a holomorphic function $F$ on the whole set $\wt{D}$ by Hartogs' Extension Theorem for singular spaces
(see \cite{MePo2}, or \cite{Rp5}), and
$$F\in \OO(\wt{D}) \subset L^\infty_{loc}(\wt{D}).$$
But then
$$f':= f- F \in L^p(D^*)$$
is the desired solution of $\dq f' =\omega$, because
$\supp f'\subset\supp\omega$ by the identity theorem
for holomorphic functions. That proves Theorem \ref{thm:compact}.


\section{Necessary Conditions (Theorem \ref{thm:necessary})}\label{sec:necessary}

We will now prove Theorem \ref{thm:necessary}. So, let $D\subset\subset Y$ be a bounded open set 
such that $0\in D$, $D^*=D\setminus\{0\}$, $G=\pi^{-1}(D)$ and $G^*=\pi^{-1}(D^*)$.
We will use the exhaustion function $\rho:=\|\cdot\|^2\circ \pi$ which is strictly plurisubharmonic on $N$ outside
the zero section $X$. Then there exist indices $\epsilon>0$ and $\delta>0$ such that
$$G_\epsilon \subset \subset G \subset \subset G_\delta,$$
where
$$G_\epsilon=\{z\in N: \rho(z)<\epsilon\}$$
and
$$G_\delta = \{z\in N:\rho(z)<\delta\}$$
are smoothly bounded strongly pseudoconvex neighborhoods of the zero section in $N$.
We can again use the fact that
$$\bigoplus_{\mu\geq c(p,q,d)} H^q(X,\OO(N^{-\mu})) \cong H^q(G_\epsilon, \mathcal{I}^{c(p,q,d)} \OO)$$
by Theorem 5.1 in \cite{Rp4}.
We must now clearify what the Definition of $c(p,q,d)$ means for us:

\begin{lem}\label{lem:c}
There exists $0<\nu<1$ such that the following is true:
Let
$$t=(q-1)-\frac{2d-2}{p} + \nu,$$
and $k(p,t)$ the $\dq$-weight of $(p,t)$. Then:
$$c(p,q,d) = k(p,t).$$
\end{lem}

\begin{proof}
When we represent the $\dq$-weight $k(p,t)$ by the formula
in Lemma \ref{lem:k1}, then
it is easy to see that there exists $0<\nu<1$ such that
\begin{eqnarray*}
k(p,t) &=& \max\left\{\begin{array}{ll}
m\in \Z: m < q + 1 -\frac{2d}{p} + \nu& , p\neq 1\\
m\in \Z: m\leq q+1 -\frac{2d}{p} + \nu&, p=1
\end{array}\right\}\\
&=& \max\{m\in\Z: m\leq q+1 -2d/p\} = c(p,q,d).
\end{eqnarray*}
One just has to choose $\nu>0$ small enough.
\end{proof}

We will abbreviate $c(p,q,d)$ by $c$ from now on.
By use of Lemma \ref{lem:c},
the exact sequence in Theorem \ref{thm:main}
tells us that a class $[\omega]\in H^q(G_\epsilon, \mathcal{I}^c \OO)$
can be represented by a form
\begin{eqnarray*}
\omega_1 \in |\mathcal{I}|^{t} \mathcal{L}^p_{0,q}(G_\epsilon).
\end{eqnarray*}
But, we will see that there also exists
\begin{eqnarray*}
\omega_2 \in |\mathcal{I}|^{t} \mathcal{L}^p_{0,q}(G_\delta).
\end{eqnarray*}
such that $[\omega_2|_{G_\epsilon}] = [\omega_1] =[\omega] \in H^q(G_\epsilon, \mathcal{I}^c\OO)$.

That follows from the following consideration:
As in the beginning of the proof of Theorem \ref{thm:sufficient}, Grauert's bump method
(see \cite{LiMi}, chapter IV.7) tells us that the mapping
\begin{eqnarray}\label{eq:bumps}
H^q(G_\delta, \mathcal{I}^c \OO) \rightarrow H^q(G_\epsilon, \mathcal{I}^c \OO)
\end{eqnarray}
induced by restriction of forms is surjective. For later use, we remark
that it is in fact an isomorphism because the groups under consideration are of equal finite dimension.\\

Now, let
\begin{eqnarray*}
s:=q -\frac{2d-2}{p} = t + (1-\nu).
\end{eqnarray*}
We will show that we can also assume
\begin{eqnarray*}
\omega_2 \in |\mathcal{I}|^s \mathcal{L}^p_{0,q}(\wt{G}),
\end{eqnarray*}
where $G\subset\subset \wt{G} \subset\subset G_\delta$.
This follows from the fact that we can solve the $\dq_c$-equation locally from
$$|\mathcal{I}|^t \mathcal{L}^p_{0,q}$$
into
$$|\mathcal{I}|^{t+1-\nu} \mathcal{L}^p_{0,q-1} = |\mathcal{I}|^s \mathcal{L}^p_{0,q-1}.$$
That is well-known at points not on the exceptional set $X=E$,
because at such points we only have to solve from $L^p_{0,q}$ into $L^p_{0,q-1}$.
At points on the exceptional set, it follows from Theorem 2.1 in \cite{Rp4},
because $c$ is the $\dq$-weight of $(p,t)$.
So, cover a domain which is slightly smaller than $G_\delta$ by finitely many domains $\{U_j\}_{j\in J}$
where we have solutions
$$\dq_c v_j = \omega_2|_{U_j}, \ \ v_j \in |\mathcal{I}|^s \mathcal{L}^p_{0,q-1}(U_j).$$
Then, let $\{\chi_j\}_{j\in J}$ be a smooth
partition of unity associated to $\{U_j\}_{j\in J}$, and define
$$\eta:=\sum_{j\in J} \chi_j v_j \in |\mathcal{I}|^s \mathcal{L}^p_{0,q-1}(\wt{G}),$$
where $\wt{G}:=\bigcup U_j$. Then, we calculate:
$$\dq_c \eta = \sum_{j\in J} \chi_j \dq_c v_j + \sum \dq \chi_j \wedge v_j= : \omega_2 - \omega.$$
Therefore $[\omega_2]=[\omega] \in H^q(G_\epsilon,\mathcal{I}^c \OO)$ can in fact be represented by
$\omega\in |\mathcal{I}|^s \mathcal{L}^p_{0,q}(\wt{G})$.
Now, it follows from the last statement of Lemma \ref{lem:blowup} that
$$(\pi|_{G\setminus X}^{-1})^* \omega \in L^p_{0,q}(D^*),$$
and it is clear that this form is $\dq$-closed in the sense of distributions on $D^*$. Hence, $\omega$ determines a class in $H^q_{(p)}(D^*,\OO)$.
We have to show that this assignment defines in fact a mapping from $H^q(G_\epsilon,\mathcal{I}^c \OO)$
into $H^q_{(p)}(D^*,\OO)$.\\

Because \eqref{eq:bumps} is an isomorphism, we only have to consider what happens
if the class $[\omega]\in H^q(G_\epsilon,\mathcal{I}^c\OO)$ is given by a different form
$$\omega_2' \in |\mathcal{I}|^t \mathcal{L}^p_{0,q}(G_\delta)$$
such that
$$\omega_2 - \omega_2' = \dq_c \theta$$
with $\theta \in |\mathcal{I}|^t \mathcal{L}^p_{0,q-1}(G_\delta)$.
Now, construct $\omega' \in |\mathcal{I}|^s\mathcal{L}^p_{0,q}(\wt{G})$ from $\omega_2'$
analogously to the construction of $\omega$ (from $\omega_2$).
Then, it follows that
\begin{eqnarray*}
\omega-\omega' &=& - \dq_c \eta + \dq_c \eta' +\omega_2 - \omega_2'\\
&=& - \dq_c \eta + \dq_c \eta' + \dq_c \theta|_{\wt{G}} = \dq_c \Delta,
\end{eqnarray*}
with
$$\Delta := \eta'-\eta + \theta|_{\wt{G}} \in |\mathcal{I}|^t \mathcal{L}^p_{0,q-1}(\wt{G}).$$
Furthermore, we get
$$(\pi|_{G\setminus X}^{-1})^* \omega - (\pi|_{G\setminus X}^{-1})^* \omega' = \dq (\pi|_{G\setminus X}^{-1})^* \Delta,$$
where $(\pi|_{G\setminus X}^{-1})^* \Delta \in L^p_{0,q-1}(D^*)$ by the last statement of Lemma \ref{lem:blowup},
because 
$$t\geq s-1=(q-1) -\frac{2d-2}{p}.$$
This shows that $\omega$ and $\omega'$ determine the same class in $H^q_{(p)}(D^*,\OO)$,
and hence our mapping is well-defined. It remains to show that it is injective.
That can be done by integration along the fibers.
So, assume that
$$(\pi|_{G\setminus X}^{-1})^* \omega = \dq \alpha$$
on $D^*$ where $\alpha \in L^p_{0,q-1}(D^*)$. Let $\chi \in C^\infty_{cpt}(D)$
be a smooth cut-off function with compact support which is identically $1$ in a neighborhood of the origin.
Then:
$$\beta := \omega - \dq \pi^* \big( (1-\chi) \alpha\big) \in |\mathcal{I}|^t \mathcal{L}^p_{0,q}(G_\epsilon)$$
has compact support in $G_\epsilon$ and is $\dq_c$-closed.
Thus, integration along the fibers of $N$ as a holomorphic line bundle over $X$ (as in section \ref{sec:fibers}) gives
$$\gamma \in |\mathcal{I}|^t \mathcal{L}^p_{0,q-1}(G_\epsilon)$$
such that
$$\dq_c\gamma=\beta,\ \ \ \mbox{ and }\ \ \dq_c \left(\gamma + \pi^* \big( (1-\chi) \alpha\big)\right)= \omega,$$
where
$$\gamma + \pi^* \big( (1-\chi) \alpha\big) \in |\mathcal{I}|^t \mathcal{L}^p_{0,q-1}(G_\epsilon).$$
Here, one should recall that $c=k(p,t)$ by Lemma \ref{lem:c}.
So, Theorem \ref{thm:main} yields $[\omega]=0 \in H^q(D_\epsilon,\mathcal{I}^c \OO)$,
as we intended to show. That completes the proof of Theorem \ref{thm:necessary}.

\section{Examples and Applications}\label{sec:examples}

As a direct consequence of Theorem \ref{thm:sufficient}, we obtain:

\begin{thm}\label{thm:bounded}
Let $X$, $Y$ and $N$ be as in the introduction, $D\subset\subset Y$ strongly pseudoconvex such that $0\in D$, $D^*=D\setminus\{0\}$,
$1\leq p\leq \infty$, $1\leq q\leq d=\dim Y$. Set
$$a(p,q,d) :=
\left\{\begin{array}{ll}
\max\{k\in\Z: k< 1+q- 2d/p\} &, p \neq 1,\\
\max\{k\in\Z: k\leq 1+ q -2d/p\} &, p=1,
\end{array}\right.$$
and assume that 
$$\bigoplus_{\mu\geq a(p,q,d)} H^q(X,\OO(N^{-\mu}))=0.$$
Then there exists a bounded linear operator
$${\bf S}_q: L^p_{0,q}(D^*)\cap \ker \dq \rightarrow L^p_{0,q-1}(D^*)$$
such that $\dq {\bf S}_q \omega= \omega$.
\end{thm}

\begin{proof}
Theorem \ref{thm:sufficient} tells us that
$$H^q_{(p)}(D^*,\OO)=0.$$
Now, the statement follows by a standard technique based on the open mapping theorem (see for example \cite{FOV1}, Lemma 4.2),
because $L^p_{0,q}(D^*)$ and $L^p_{0,q-1}(D^*)$ are Banach spaces.
\end{proof}

Let us take a look at two simple examples in the case $d=\dim Y=2$.
Then, $X$ is a compact Riemann surface.
Firstly, let us assume that genus $g(X)=0$, hence $X\cong \C\mathbb{P}^1$.
Let $z_0\in X$ be an arbitrary point 
and $D= - (z_0)$ the associated divisor.
Then it follows that
\begin{eqnarray*}
H^j(X,\OO(\mu D))\cong H^j(X,\OO(N^\mu))
\end{eqnarray*}
for all $j\geq 0$ and $\mu\in\Z$.
It is well-known (and easy to calculate by power series)
that
$$l(\mu):=\dim H^0(X,\OO(-\mu D)) = \left\{
\begin{array}{ll}
1 + \mu&,\mbox{ for } \mu\geq -1,\\
0&,\mbox{ for } \mu \leq -1,
\end{array}\right.$$
because $H^0(X,\OO(-\mu D))$ is the space of meromorphic functions on $X$ with a single pole of order $\mu$ at $z_0$.
Hence, we calculate by the Theorem of Riemann-Roch that
\begin{eqnarray*}
- \dim H^1(X,\OO(N^{-\mu})) = \deg(-\mu D) + 1 - g(\C\mathbb{P}^1) - l(\mu) =
\left\{
\begin{array}{ll}
0&,\mbox{ for } \mu\geq -1,\\
-1&,\mbox{ for } \mu = -2.
\end{array}\right.
\end{eqnarray*}
Therefore, Theorem \ref{thm:sufficient} and Theorem \ref{thm:necessary} tell us that
\begin{eqnarray*}
\dim H^1_{(p)} (D^*,\OO) \left\{
\begin{array}{ll}
=0&, \mbox{ if } p>4/3,\\
\leq 1&, \mbox{ if } p=4/3,\\
=1&, \mbox{ if } p<4/3,
\end{array}\right.
\end{eqnarray*}
if $D$ is a strongly pseudoconvex neighborhood of the origin in $Y$.
An important example for such a variety is $Y=\{(x,y,z)\in\C^3: xy=z^2\}$.\\

As a second example, we use the same construction, but assume that $X$ is an elliptic curve in $\C\mathbb{P}^{n-1}$.
Here, $H^0(X,\OO(-\mu D))$ is the space of elliptic functions with a single pole of order $\mu$.
So, it is well-known that we have
$$l(\mu):=\dim H^0(X,\OO(-\mu D)) =\left\{
\begin{array}{ll}
0 &, \mbox{ for } \mu \leq -1,\\
1 &, \mbox{ for } \mu = 0,\\
\mu &, \mbox{ for } \mu\geq 1.
\end{array}\right.$$
Using the Theorem of Riemann-Roch again, we calculate
\begin{eqnarray*}
- \dim H^1(X,\OO(N^{-\mu})) = \deg(-\mu D) + 1 - g(X) - l(\mu)
= \left\{
\begin{array}{ll}
\mu &, \mbox{ for } \mu \leq -1,\\
-1 &, \mbox{ for } \mu = 0,\\
0 &, \mbox{ for } \mu\geq 1.
\end{array}\right.
\end{eqnarray*}

Therefore, Theorem \ref{thm:sufficient} and Theorem \ref{thm:necessary} tell us that
\begin{eqnarray*}
\dim H^1_{(p)} (D^*,\OO) \left\{
\begin{array}{ll}
=0&, \mbox{ if } p>4,\\
\in\{0, 1\}&, \mbox{ if } p=4,\\
=1&, \mbox{ if } 4>p>2,\\
\in \{1,2\}&, \mbox{ if } p=2,\\
=2&, \mbox{ if } 2>p> 4/3,\\
\in\{2,3,4\} &, \mbox{ if } p=4/3,\\
=4&, \mbox{ if } 4/3> p,
\end{array}\right.
\end{eqnarray*}
if $D$ is a strongly pseudoconvex neighborhood of the origin.
Examples are the varieties $Y=\{(x,y,z)\in\C^3: y^2z=x^3+axz^2 + bz^3\}$
for suitable values of $a$, $b$.\\

Let us return to the first example $X \cong \C\mathbb{P}^1$.
Combining that consideration with Theorem \ref{thm:necessary}, we obtain:

\begin{thm}\label{thm:l2}
Let $X$ and $Y$ be as in the introduction, $\dim Y=2$ and
$D\subset\subset Y$ strongly pseudoconvex such that $0\in D$, $D^*=D\setminus\{0\}$.
Then:
$$H^1_{(2)}(D^*,\OO)=0\ \ \Leftrightarrow\ \ X\cong \C\mathbb{P}^1.$$
\end{thm}

\begin{proof}
By assumption, $X$ is a compact Riemann surface. If $X\cong \C\mathbb{P}^1$ then $H^1_{(2)}(D^*,\OO)=0$
by the considerations above. Conversely, if $H^1_{(2)}(D^*,\OO)=0$, then $H^1(X,\OO)=0$ 
by Theorem \ref{thm:necessary}, and that implies that $X\cong\C\mathbb{P}^1$.
\end{proof}

This example, namely the groups $H^1_{(2)}$ at isolated singularities of co-dimension two, are of special interest
because of the following Extension Theorem of Scheja (see \cite{Sch1,Sch2}), which settles the case of higher co-dimension:

\begin{thm}\label{thm:scheja1}
Let $Y$ be a closed pure dimensional analytic subset in $\C^n$ which is locally a complete
intersection, and $A$ a closed pure dimensional analytic subset of $Y$.
Then, the natural restriction mapping
$$H^q(Y,\OO_Y) \rightarrow H^q(Y\setminus A,\OO_{Y\setminus A})$$
is bijective for all $0\leq q\leq \dim Y-\dim A-2$.
\end{thm}

Using this result, our integration along the fibers yields:

\begin{thm}\label{thm:scheja2}
Let $X$ and $Y$ be as in the introduction, and $D\subset\subset Y$ strongly pseudoconvex such that $0\in D$, $D^*=D\setminus\{0\}$,
$1\leq p\leq \infty$, and
$$1\leq q\leq d-2=\dim Y-2.$$
Then:
$$H^q_{(p)}(D^*,\OO)=0.$$
\end{thm}

\begin{proof}
As in the beginning of section \ref{sec:sufficient},
assume that 
$$D\cap U =\{z\in U:\rho(z)<0\}$$ 
where $\rho\in C^2(U)$ is a regular strictly plurisubharmonic 
defining function on a neighborhood $U$ of $bD$, and that
there exists $\epsilon>0$ such that 
$$D_\epsilon:= D\cup \{z\in U:\rho(z)<\epsilon\}$$ 
is a strongly pseudoconvex extension of $D$.
So, it follows by Grauert's bump method that the natural homomorphism
\begin{eqnarray}\label{eq:bump11}
r_q: H^q_{(p)}(D_\epsilon^*,\OO) \rightarrow H^q_{(p)}(D^*,\OO)
\end{eqnarray}
(induced by restriction of forms) is surjective (see \cite{LiMi}, chapter IV.7).
Here, we also set $D_\epsilon^*=D_\epsilon\setminus\{0\}$.
We need to observe that $D_\epsilon$ is a Stein domain. But that follows from the fact 
that $D_\epsilon$ is a bounded strongly pseudoconvex domain in the Stein space $Y$ (see \cite{Na2}).
Moreover, $Y$ is a complete intersection, and so Theorem \ref{thm:scheja1} tells us that
\begin{eqnarray}\label{eq:stein}
H^q(D^*_\epsilon,\OO)=0.
\end{eqnarray}
Now, let $[\omega]\in H^q_{(p)}(D^*,\OO)$ be represented by the $\dq$-closed form $\omega\in L^p_{0,q}(D^*)$.
Because \eqref{eq:bump11} is surjective, we can assume that
$$\omega\in L^p_{0,q}(D^*_\epsilon).$$
But now \eqref{eq:stein} tells us that there exists
$$\eta\in L^p_{(0,q-1),loc}(D^*_\epsilon)$$
such that $\dq \eta=\omega$ on $D^*_\epsilon$.
So, choose a smooth cut-off function $\chi\in C^\infty_{cpt}(D)$ with compact support in $D$ such that
$\chi$ is identically $1$ in a neighborhood of the origin.
It follows that $[\omega]\in H^q_{(p)}(D^*,\OO)$ can be represented by
$$\tau:= \omega - \dq\big((1-\chi) \eta\big) \in L^p_{0,q}(D^*)$$
which has compact support in $D$. The integration along the fibers which we have already used
in the sections \ref{sec:fibers} and \ref{sec:necessary} tells us that we can use the
blow up of $Y$ to produce a solution
$$\sigma \in L^p_{0,q-1}(D^*)$$
such that $\dq \sigma=\tau$ (see Theorem \ref{thm:integration}), and that finishes the proof.
\end{proof}


Combining Theorem \ref{thm:scheja2} with Theorem \ref{thm:necessary} (for $p=1$), we obtain immediately:

\begin{thm}\label{thm:scheja3}
Let $X$, $Y$ and $N$ be as in the introduction, and
$$1\leq q\leq d-2=\dim Y-2.$$ 
Then it follows that
$$H^q(X,\OO(N^{-\mu}))=0$$
for all $\mu\geq 1+q -2d$.
\end{thm}

Similarly, we can deduce from Theorem \ref{thm:necessary}:

\begin{thm}\label{thm:extension}
Let $N$ be the universal bundle over $\C\mathbb{P}^{k}$ for $k\geq 1$,
and $1\leq q \leq k$. Then:
$$H^q(\C\mathbb{P}^k, \OO(N^{-\mu}))=0\ \ \ \mbox{ for all } \mu\geq q-2k.$$
\end{thm}

\begin{proof}
Note that we have at no place assumed that $X$ is a proper subset of $\C\mathbb{P}^{n-1}$,
respectively that $Y$ should be a proper subset of $\C^n$.
So, in the setting of the introduction let
$X=\C\mathbb{P}^{n-1}=\C\mathbb{P}^k$
and
$Y=\C^n=\C^{k+1}$.
Moreover, let $D$ be the unit ball in $Y$ and $p=2n/(2n-1)$. Now, if 
$\omega \in L^p_{0,q}(D^*)$
is a $\dq$-closed form on the punctured ball, then the $\dq$-Extension Theorem 3.2 in \cite{Rp4}
tells us that in fact $\omega$ defines a $\dq$-closed $L^p$-form on the whole ball.
So, there exists $\eta\in L^p_{0,q-1}(D)$ such that $\dq\eta=\omega$ (see \cite{Kr}), and it follows that
$H^q_{(p)}(D^*,\OO)=0$
for all $1\leq q \leq k+1$. Thus, Theorem \ref{thm:necessary} implies that
$H^q(\C\mathbb{P}^k, \OO(N^{-\mu}))=0$
for all
\begin{eqnarray*}
\mu &\geq& q+1 -\frac{2n}{p} = q+1 - 2n + 1 = q-2k.
\end{eqnarray*}
\end{proof}

{\bf Acknowledgments}

\vspace{2mm}
This work was done while the author was visiting the University of Michigan at Ann Arbor,
supported by a fellowship within the Postdoc-Programme of the German Academic Exchange Service (DAAD).
The author would like to thank the SCV group at the University of Michigan for its hospitality.


\end{document}